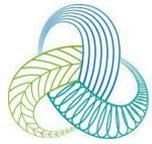 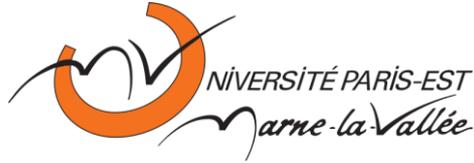 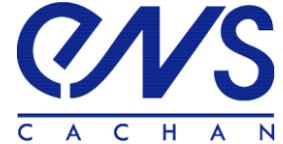

# The Logit lane assignment model: first results.

Nadir Farhi[a],
Habib Haj-Salem[a],
Megan Khoshyaran[b],
Jean-Patrick Lebacque[a,*],
Francesco Salvarani[c],
Bernard Schnetzler[a],
Florian de Vuyst[c]

[a] Université Paris-Est, IFSTTAR, GRETTIA, F-93166 Noisy-le-Grand, France.
[b] ETC Economics Traffic clinic, France.
[c] ENS Cachan (Ecole Normale Supérieure de Cachan), France.

* `jean-patrick.lebacque@ifsttar.fr` (Corresponding author).

N. Farhi, H. Haj-Salem, M.M. Khoshyaran, J.P. Lebacque, F. Salvarani, B. Schnetzler, F. de Vuyst

**ABSTRACT**

The Logit lane assignment model has been introduced recently in order to describe multi-lane traffic flow from a macroscopic point of view. The model is based on the idea that each available lane has a specific utility for each driver, who chooses the lane with the highest utility. The model is expressed by a system of conservation laws with a smooth but implicitly defined flux function. The first aim of the paper is to explore on two data-sets how traffic data supports the fact that traffic speed constitutes an explanatory variable of lane assignment. Second the paper addresses the problem of discretization of the model. Several numerical schemes are proposed: Lax-Friedrichs, Euler-lagrange remap, Lagrange, and their convergence properties are illustrated on the treatment of the Riemann problem. Directions for future research are outlined

N. Farhi, H. Haj-Salem, M.M. Khoshyaran, J.P. Lebacque, F. Salvarani, B. Schnetzler, F. de Vuyst

2   **INTRODUCTION**

3   During the last decade, in order to enhance the mobility in the transportation network, several
4   control actions are implemented and tested aiming at improving traffic condition of the
5   overall network. In particular, modern transportation management systems are based more
6   and more on the application of the new technologies of Intelligent Transport System (ITS)
7   Therefore contemporary traffic control measures aim to pursue a more efficient use of
8   existing infrastructure by recognizing the importance of distinguishing user-classes, such as
9   trucks, buses, passenger-cars, and high-occupancy vehicles, and by managing these classes
10  differently.  Among the control actions the managed-lanes take an increasing importance (see
11  for instance Cohen 2001). High occupancy vehicle (HOV) lanes and high occupancy toll
12  (HOT) lanes are the two most common managed lanes which people see daily. However,
13  before implementing these lane assignment control actions, off-line simulation studies are
14  performed aiming at evaluating the traffic impact before the implementation phase.  In this
15  case, the general used traffic models are of microscopic type. This permits to simulate the
16  behavior of the driver, lane changing etc. with more or less accuracy. A very large literature
17  on this topic is available.

18  This paper focuses on the modeling of multilane macroscopic traffic flow. Most of these
19  models are restricted to uni-directional freeway traffic and treat the different lanes of a road in
20  terms of global quantities, such as the capacity and the possibilities for overtaking. However,
21  this kind of simplification is clearly not applicable if there is disequilibrium between
22  neighboring lanes. Therefore, some researchers carried out empirical investigations of the
23  observed density oscillations between neighboring lanes or proposed models for their mutual
24  influences. However, these models are of phenomenological type and treat the inter-lane
25  interactions in a rather heuristic way. Moreover, most of them are based on the simple traffic
26  flow model of Lighthill, Whitham and Richards (LWR, Lighthill and Whitham 1955,
27  Richards 1956) which assumes the average velocity on each lane to be in equilibrium with the
28  density. The simplest models consider the global impact of lane change without actually
29  modeling the lanes (Greenberg et al 2003, Jin 2006 for instance). We also mention
30  Michalopoulos et al. 1984 which is based on a relaxation process tending to equilibrate lane
31  densities. An approach describing the flows between lanes (governed by differences in speed)
32  and the impacts of lane change, has been developed in a series of papers, (we refer for
33  instance to Laval and Daganzo 2005 and 2006, Laval and Leclercq 2008, see also Schnetzler
34  et al 2010). Kinetic models propose an approach which integrates microscopic lane change
35  behavior into macroscopic models; which has been introduced by Helbing 1997, and has been
36  developed by Hoogendoorn and Bovy 1999, Ngoduy 2006 with careful inclusion of lane
37  change motives.

38  At an intermediate level a relatively precise description can be achieved by considering the
39  net result of flows between lanes without describing the details of the dynamics which lead to
40  this net result. If we consider only the traffic state on lanes, we can assume that there is
41  equilibrium between lanes. Daganzo 1997 initiated this approach in the case of special lanes
42  and two classes of vehicles based on an all or nothing assignment. Lebacque and Khoshyaran
43  1998-2002 proposed an extension to many lanes and vehicle classes. This model is expressed

N. Farhi, H. Haj-Salem, M.M. Khoshyaran, J.P. Lebacque, F. Salvarani, B. Schnetzler, F. de Vuyst

1  as a system of conservation laws, the flux function of which is continuous but only defined
2  implicitly, and admits only a piece-wise continuous gradient.

3  The model presented in this paper generalizes the ideas of Lebacque and Khoshyaran 2009. It
4  is based on the idea that at a macroscopic scale, it is possible to neglect the microscopic
5  details of the lane-change maneuvers, and to model only the net result of these maneuvers.
6  Further, the model assumes that this result can be described as a user equilibrium in which the
7  lanes constitute the alternatives. The model assumes a stochastic utility describing the lane
8  choice process. This model is both more regular and tractable, and better supported by data.

9  This paper is organized as follows: section 1 is dedicated to a description of the essential
10 features of the used model. Section 2 presents some results of data measurements analyzed in
11 the perspective of the model. Section 3 describes three numerical schemes, the classical Lax-
12 Friedrichs scheme, the Euler-remap scheme and a Lagrangian scheme.

13 **THE MODEL**

14 We are interested here by modeling the traffic on a multilane road. Cars move in one direction
15 on a multilane road. Lane change is allowed, and passing is permitted by lane change. Drivers
16 are classified in several classes according to their destination. This assumption is important
17 because the lane assignment of drivers may depend only on their classes. For example,
18 approaching a divergent, some drivers have to change lane in order to take the desired
19 destination, independent of the car density on each lane.

20 We use the following notations:

21 - $I$ : the set of lanes indexed by $i$.
22 - $D$ : the set of driver classes indexed by $d$.
23 - $I^d$ : the set of lanes accessible to user class $d$.
24 - $D_i$ : the set of user classes that can move on lane $i$.
25 - $\rho^d(x,t)$ : the car-density of class $d$ users at location $x$ and time $t$.
26 - $\rho_i(x,t)$ : the car-density at location $x$ on lane $i$ at time $t$.
27 - $\rho_i^d(x,t)$ : the car density of class $d$ users at location $x$ on lane $i$ at time $t$.
28

29 Note that with these notations $\rho^d = \sum_{i \in I^d} \rho_i^d$ and $\rho_i = \sum_{d \in D_i} \rho_i^d$.

30 The car density $\rho(x,t)$ at location $x$ and time $t$ is then given by

$$\rho = \sum_{d \in D} \rho^d = \sum_{d \in D} \sum_{i \in I^d} \rho_i^d = \sum_{i \in I} \sum_{d \in D_i} \rho_i^d = \sum_{i \in I} \rho_i.$$

31 Utilities are associated to the lanes of the road, depending on the speed of traffic on each lane.
32 Utilities are user-class specific. The speed of traffic on a lane depends on the car-density on
33 the considered lane, according to the fundamental traffic diagrams of each lane. The choice of
34 lanes by drivers is based on the utility of each lane. However, since the utility of a lane
35 decreases when the density increases, a massive choice of a lane with a high utility has the net
36 effect of lowering the car speed on that lane, inducing an equilibrium. The model is

N. Farhi, H. Haj-Salem, M.M. Khoshyaran, J.P. Lebacque, F. Salvarani, B. Schnetzler, F. de Vuyst

macroscopic: therefore we neglect the actual process of lane change and only consider the macroscopic result of this process. We assume here that the traffic dynamics on the road consists in splitting, at each time step and at each location, the car density over all the lanes by satisfying a user equilibrium. We notice that the car-densities are conservative by user class, but not by lane.

The utility $U_i^d$ of a lane $i \in I$ for a user of class $d \in D$ is

$$U_i^d = v_i + \theta_i^d + \xi_i^d,$$

where $v_i = V_i(\rho_i)$ is the car-speed on lane $i$, given by the fundamental traffic diagram $V_i$ on that lane, $\theta_i^d$ is a constant expressing the preference of class $d$ users for lane $i$, and $\xi_i^d$ is a Gumbel random variable that expresses the stochastic elements of the lane choice. Users of each class chose the lane with the highest utility.

The car-density is split on the different lanes according to the Logit lane assignment model

$$\rho_i^d = \rho^d \frac{\exp\left(\{V_i(\sum_{d \in D_i} \rho_i^d) + \theta_i^d\}/\nu\right)}{\sum_{j \in I^d} \exp\left(\{V_j(\sum_{d \in D_j} \rho_j^d) + \theta_j^d\}/\nu\right)}, \quad \forall d \in D, \forall i \in I^d \quad (1)$$

with $\rho^d = \sum_{j \in I^d} \rho_j^d \quad \forall d \in D,$ and $\nu$ a sensitivity parameter.

Then at every time step a fixed point problem is resolved with respect to the $\rho_i^d$ in order to express the $\rho_i^d$ as functions of the partial densities $\rho^d$; see (Khoshyaran & Lebacque, 2012) for more details.

The partial densities $\rho_i^d$ are solutions of the following linear concave optimization problem:

$$\max_{(\rho_i^d)_{d \in D, i \in I^d}} \sum_{i \in I^d} \int_0^{\rho_i} V_i(r) dr + \sum_{d \in D, i \in I^d} \rho_i^d \mu_i^d - \nu \sum_{d \in D, i \in I^d} \rho^d H\left(\frac{\rho_i^d}{\rho^d}\right)$$
$$\left| \rho^d = \sum_{i \in I^d} \rho_i^d \quad \forall d \in D \right. \quad (2)$$

where $H(x) \stackrel{\text{def}}{=} x[\ln x - 1]$ is the negentropy function. This can be checked by stating the lagrange optimality conditions for (1). Another way to explain this result is to recall that the Beckmann transformation recasts any assignment problem with increasing costs as an optimization problem.

It follows that the partial densities $\rho_i^d$ and the densities per lane $\rho_i$ can be expressed as implicit but smooth (indefinitely differentiable) functions of the densities per destination $\bar{\rho} \stackrel{\text{def}}{=} (\rho^d)_{d \in D}$:

$$\rho_i^d \stackrel{\text{def}}{=} \mathfrak{R}_i^d(\bar{\rho}), \quad \rho_i \stackrel{\text{def}}{=} \sum_{\delta \in D_i} \rho_i^\delta \stackrel{\text{def}}{=} \mathfrak{R}_i(\bar{\rho}), \quad (3)$$

and that the flow $q^d$ of densities per destination $d$ is given by:

$$q^d \stackrel{\text{def}}{=} \mathfrak{H}^d(\bar{\rho}) = \sum_{i \in I_d} \rho_i^d V_i(\rho_i) = \sum_{i \in I_d} \mathfrak{R}_i^d(\bar{\rho}) V_i\left(\mathfrak{R}_i(\bar{\rho})\right). \quad (4)$$

N. Farhi, H. Haj-Salem, M.M. Khoshyaran, J.P. Lebacque, F. Salvarani, B. Schnetzler, F. de Vuyst

1  The flux function $\overline{\mathfrak{H}} \stackrel{\text{def}}{=} (\mathfrak{H}^d)_{d \in D}$ is smooth but is defined implicitly. It is possible to calculate
2  the gradient and eigenvalues of the flux function, but the calculations are quite technical and
3  beyond the scope of this paper.

4  The Logit lane assignment traffic flow model can now be expressed as the following system
5  of conservation laws:

6  $$\partial_t \bar{\rho} + \partial_x \overline{\mathfrak{H}}(\bar{\rho}) = 0 \text{ or equivalently } \partial_t \rho^d + \partial_x \mathfrak{H}^d(\bar{\rho}) = 0 \quad \forall d \in D. \quad (5)$$

7  The reasons for which this system is difficult to solve are, in increasing order of difficulty:
8  - The implicit definition of the flux functions,
9  - The impossibility of estimating the gradient of the flux function explicitly,
10 - The impossibility of estimating the eigenvalues of the flux function other than
11   numerically.
12

13 **TESTING THE LANE CHOICE MODEL ON DATA**

14 The model presented in the previous section does not recapture all the complexity of
15 lane change at the microscopic level. The model aims at reproducing the net result of the lane
16 assignment process at a macroscopic level. Thus validating the model on real data is difficult.
17 Essentially what is required is a comparison between measurements and model predictions.
18 As a first step, we have checked whether the lane speed is indeed a relevant explanatory
19 variable of lane assignment and whether the Logit assignment is compatible with
20 measurements. This is the object of the present section.

21 **The MARIUS data**.

22 The objective here is to determine from observations in what proportions a lane choice is
23 explained by the speed or by the random component of the utility (i.e. by unknown
24 explanatory variables).

25 The study is conducted for various levels of density, in order to assess the importance of this
26 level.

27 This determination requires the comparison between the experimental distribution of traffic
28 $p_i = \rho_i / \rho$, where $\rho_i$ is the density on the lane $i$ and where $\rho = \sum \rho_i$, and the theoretical
29 distribution given by the logit model $p_i = \exp(\nu^{-1} v_i) / \sum_k \exp(\nu^{-1} v_k)$ where $v_k$ is the observed
30 speed on the lane $k$. This computation requires only the observed speeds and it is therefore not
31 necessary to know the fundamental diagrams per lane $V_i(\rho_i)$.

32 Essentially, this approach will be valid in the homogeneous case, i.e. we can neglect the
33 preference parameters $\theta_i^d$ and assume that $I_d = I$ for all $d$. In other words: we neglect the
34 impact of on- and off-ramps on the lane choice. It is recalled that in the homogeneous case the
35 logit model is based on the assumption of a utility $U_i = v_i + \xi_i$ where $v_i$ is the known part of
36 the utility (which is identified with the speed) and $\xi_i$ is the contribution of the unknown
37 explanatory variables. The parameter $\nu$ is proportional to the standard-deviation of the
38 distributions of variables $\xi_i$. When the standard-deviation tends to infinity, the speed does not
39 explain the lane choice; vehicles are distributed with equal probability on all lanes (i.e.
40 densities are equal). So the experimental distribution is compared to the theoretical
41 distribution for values of $\nu$ ranging between 17 km/h and 2500 km/h. The best value

N. Farhi, H. Haj-Salem, M.M. Khoshyaran, J.P. Lebacque, F. Salvarani, B. Schnetzler, F. de Vuyst

minimizes the sum of the squares of the differences between the observed and the theoretical probabilities (see figure 1).

This method has been implemented on suburban highways around the city of Marseille. The selected section (A7-TC) has three lanes and is 9 km in length, with few entries or exits. It is equipped with 17 measuring stations. Each station has one sensor per lane. Each vehicle passage generates an event recording the lane, the moment of passage, the speed and the length of the vehicle. Given these event data, it is possible to reconstruct flows, densities and average speeds for different aggregation levels. Two set of daily data with respective aggregation level of 1 and 6 minutes (see table 1 and 2) have been generated. Sample sizes are large (of the order of 1000). Estimations are similar whatever the aggregation level is. Finally data sets have been recorded for 30 days. Since the estimations have no seasonality and are very similar, for any day of the recording period, the presented results are those obtained on the monthly data file.

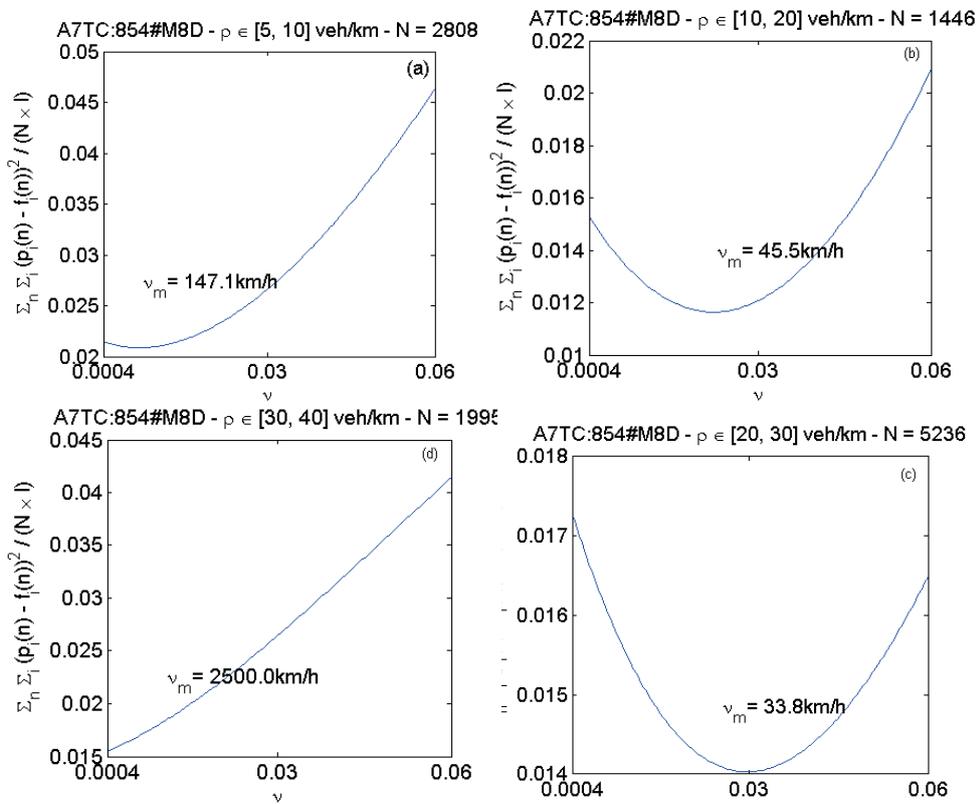

FIGURE 1. Estimates of $\nu$ on station D. Horizontal axis: $\nu$, vertical axis: normalized quadratic error

Given these preliminary results, the speed is not an explanatory variable when the traffic is fully free-flow or congested. However it is an explanatory variable for intermediate densities. Further for intermediate densities, the stochastic assignment model fits the data better than the deterministic one (which corresponds to $\nu = 0$). The values of $\nu$ are widely scattered. This scatter is not characteristic of the behavioral variability but of the difficulty of the parameter estimation. Indeed the method is not quite satisfactory (the lane choices are individual choices and not mean choices) and the sample sizes are too small in view of the complexity of the lane change interactions.

N. Farhi, H. Haj-Salem, M.M. Khoshyaran, J.P. Lebacque, F. Salvarani, B. Schnetzler, F. de Vuyst

Further, the data is at best sampled every minute, and it is not known whether the homogeneity assumption is satisfied correctly.

Assuming a probability $p_i = \exp(\nu^{-1} v_i) / \sum_k \exp(\nu^{-1} \langle v_k \rangle)$, where $v_i$ is an observed vehicle velocity and $\langle v_k \rangle$ is a mean lane velocity, a logit regression could provide better estimates. Moreover the model depends on an unsatisfied assumption: the traffic is not homogeneous. Particularly, there are lane preferences depending on the vehicle class such as truck and personal car.

To be precise, if the Logit model is satisfied and the traffic flow is heterogeneous, we would deduce from (1) by some straightforward algebra:

$$\nu \ln \frac{\rho_i}{\hat{\rho}_{j,i}} = v_i - v_j$$

with $\hat{\rho}_{j,i} = \sum_{\delta \epsilon D_i} e^{(\theta_i^\delta - \theta_j^\delta)/\nu}$, instead of $\nu \ln \frac{\rho_i}{\rho_j} = v_i - v_j$. Of course, in the homogeneous case all $\theta_i^\delta$ are equal and $\hat{\rho}_{j,i} = \rho_j$.

Finally, let us stress the fact that the utility function should depend on the speed attribute, but this dependency need not be linear.

**TABLE 1.** $\nu$ **values (km/h) for successive stations (1 mn data).**

| Densities (veh/km) | STATIONS | | | | | | | | | | | | | | | | |
|---|---|---|---|---|---|---|---|---|---|---|---|---|---|---|---|---|---|
| | A | B | C | D | E | F | G | H | I | J | K | L | M | N | O | P | Q |
| 5-10 | ∞ | ∞ | ∞ | 147 | ∞ | ∞ | ∞ | ∞ | ∞ | ∞ | ∞ | ∞ | ∞ | ∞ | 40 | 52 | 36 |
| 10-20 | ∞ | 96 | ∞ | 45 | ∞ | ∞ | 125 | 227 | 147 | ∞ | 81 | 1250 | ∞ | ∞ | 36 | 31 | 104 |
| 20-30 | ∞ | 30 | 60 | 34 | 78 | 56 | 30 | 34 | 35 | 70 | 227 | 48 | ∞ | ∞ | 33 | ? | ∞ |
| 30-40 | ∞ | ∞ | ∞ | ∞ | 86 | 57 | ∞ | ∞ | 500 | ∞ | ∞ | ∞ | ∞ | ∞ | 24 | 18 | ∞ |

N. Farhi, H. Haj-Salem, M.M. Khoshyaran, J.P. Lebacque, F. Salvarani, B. Schnetzler, F. de Vuyst

2   **TABLE 2.** $v$ **values (km/h) for successive stations (6 mn data).**

| Densities (veh/km) | STATIONS | | | | | | | | | | | | | | | | |
|---|---|---|---|---|---|---|---|---|---|---|---|---|---|---|---|---|---|
| | A | B | C | D | E | F | G | H | I | J | K | L | M | N | O | P | Q |
| 5-10 | ∞ | ∞ | ∞ | 179 | ∞ | ∞ | ∞ | ∞ | ∞ | ∞ | ∞ | ∞ | ∞ | ∞ | 39 | 51 | 34 |
| 10-20 | ∞ | 70 | 416 | 39 | ∞ | 833 | 104 | 156 | 125 | ∞ | 71 | 625 | ∞ | ∞ | 34 | 29 | 104 |
| 20-30 | ∞ | 19 | 34 | 28 | 62 | 41 | 22 | 25 | 30 | 50 | 179 | 35 | 192 | ∞ | 29 | ? | ∞ |
| 30-40 | ∞ | ∞ | ∞ | ∞ | 78 | 30 | ∞ | ∞ | ∞ | ∞ | ∞ | ∞ | ∞ | 139 | ? | ∞ | ∞ |

3
4
5   **The Boulevard périphérique data**
6
7   The object of this study is to check the Logit form of the model.
8   The site is depicted by Figure 2. The data collected comprises flows and occupancies per lane
9   every 20 seconds. Stations are available every 500 meters approximately. On- and off-ramps
10  are numerous (every 500 – 1500 meters). Thus the traffic cannot be considered homogeneous.

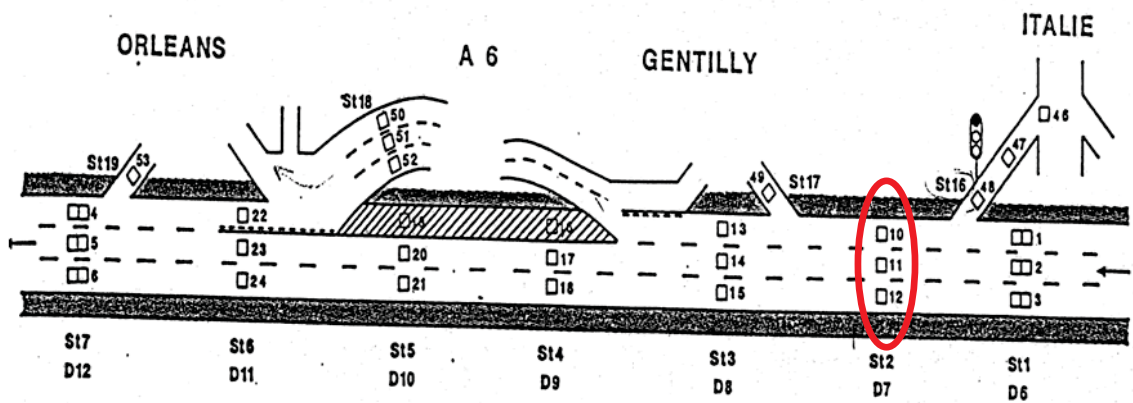

12                          FIGURE 2: Boulevard Périphérique site

N. Farhi, H. Haj-Salem, M.M. Khoshyaran, J.P. Lebacque, F. Salvarani, B. Schnetzler, F. de Vuyst

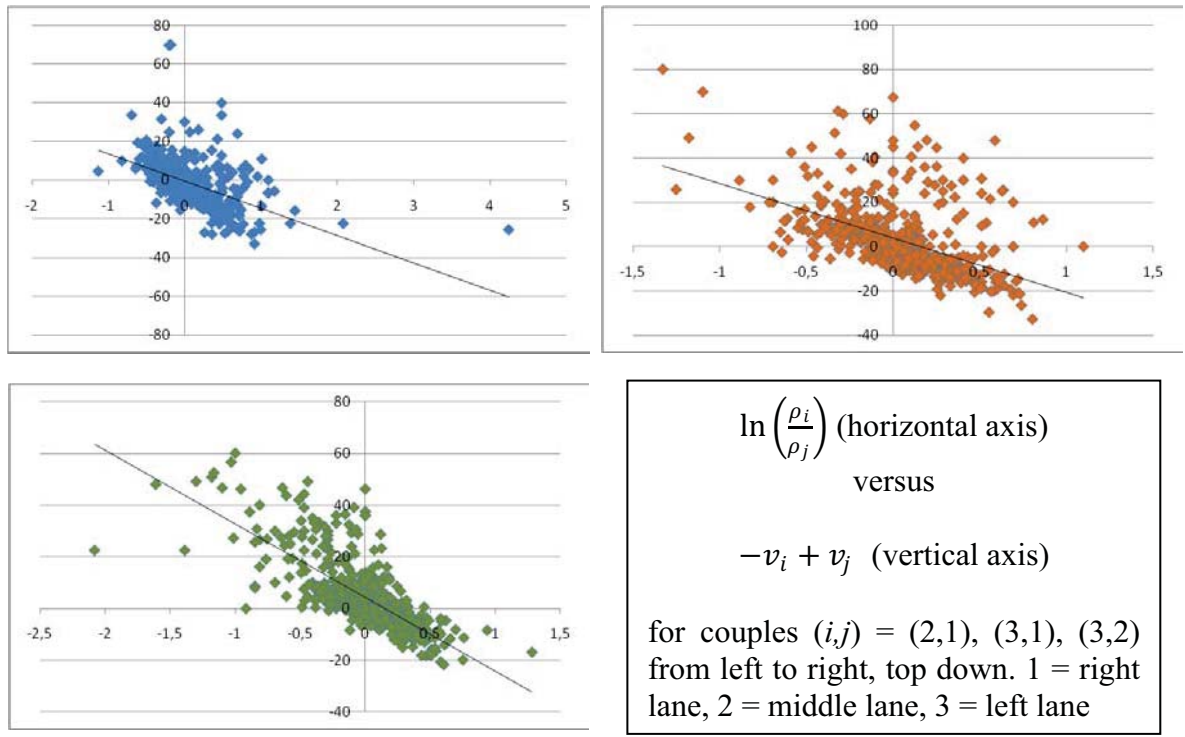

FIGURE 3: lane densities versus lane speeds. Estimation of $\nu$

In this study we have compared the values of $\ln\left(\frac{\rho_i}{\rho_j}\right)$ (for $i, j$ lanes) at various detector stations with the differences in lane speeds $v_i - v_j$. The Figure 3 shows the resulting plots and the linear regressions for the station S12 (detectors 10, 11, 12). All levels of density are included in a plot. In spite of the fact that the traffic cannot be considered homogeneous, and that the data comprises density levels for which the lane velocity has low explanatory power as shown on the Marius data, the regressions yield relatively consistent results for $\nu$. For instance the data of station S12 of Figure 3 yields the following values for $\nu$: 15, 22 and 25 km/hr. Thus the Boulevard Périphérique data supports the Logit lane assignment within the limits of the study (inhomogeneous traffic).

**NUMERICAL SCHEMES.**

The aim of this section is the development of numerical schemes for the lane assignment model. Numerical schemes should make it possible to validate the lane assignment model in a general setting (inhomogeneous traffic) and to use the model for traffic management.

**The Lax-Friedrichs scheme.**

This classic scheme will serve as a reference. The scheme is known to be first order convergent under mild conditions.

The notations are the following.

N. Farhi, H. Haj-Salem, M.M. Khoshyaran, J.P. Lebacque, F. Salvarani, B. Schnetzler, F. de Vuyst

- $\Delta t$: the time-step,
- $\Delta x$: the space discretization step,
- $x_\ell$ : the mesh (fixed) with step-size $\Delta x$,
- $\lambda \stackrel{\text{def}}{=} \Delta t / \Delta x$ : a coefficient taken to be less than the greatest eigenvalue of $\nabla \overline{\mathfrak{H}}$ (to be evaluated numerically),
- $\rho_\ell^{d,t}$ : the estimated value of the density per destination $\rho^d$ at the mesh point $x_\ell$ at time t $\Delta t$.
- $\bar{\rho}_\ell^t = \left(\rho_\ell^{d,t}\right)_{d \in D}$ : the vector of densities at location $x_\ell$ .

The Lax-Friedrichs scheme can be expressed as:

$$\rho_\ell^{d,t+1} = \tfrac{1}{2}\left(\rho_{\ell-1}^{d,t} + \rho_{\ell+1}^{d,t}\right) - \tfrac{\lambda}{2}\left(\mathfrak{H}^d(\bar{\rho}_{\ell+1}^t) - \mathfrak{H}^d(\bar{\rho}_{\ell-1}^t)\right) \qquad (6)$$

or equivalently as:

$$\bar{\rho}_\ell^{t+1} = \tfrac{1}{2}(\bar{\rho}_{\ell-1}^t + \bar{\rho}_{\ell+1}^t) - \tfrac{\lambda}{2}\left(\overline{\mathfrak{H}}(\bar{\rho}_{\ell+1}^t) - \overline{\mathfrak{H}}(\bar{\rho}_{\ell-1}^t)\right) \qquad (7)$$

The reader can refer to Holden-Risebro 2002 or Godlewski-Raviart 1996 for general results on numerical schemes for systems of conservation laws.

The quantities $\mathfrak{H}^d$ (i.e. the fluxes) are calculated by solving (1) or (2) with respect to the partial densities $\mathfrak{R}_i^d(\bar{\rho}_{\ell-1}^t)$ and $\mathfrak{R}_i^d(\bar{\rho}_{\ell+1}^t)$, and then by applying (3) and (4). In the case of two destinations and two groups of lanes only, a simple Newton iteration suffices to solve the fixed point (1). In the case of a higher dimensional problem, a projected Newton optimization algorithm can be applied to the smooth concave linear optimization problem (2) in order to estimate the partial densities.

It must be emphasized here that the evaluation of the fluxes $\mathfrak{H}^d$ constitutes a step with a much higher computational cost than the cost of all other steps of the algorithm. On the other hand the implicit evaluation of the fluxes does not affect the convergence because this step and can be performed with arbitrary precision.

**Euler-lagrange remap scheme.**

This scheme is based on a piecewise constant approximation of the solution. The notations for this scheme are the following.
- $x_\ell$ : the mesh (fixed) with step-size $\Delta x \stackrel{\text{def}}{=} h$,
- $\Delta t$ : the time-step,
- $v_\ell^{d,t}$: the velocity of the vehicule $N_\ell^{d,t}$ with destination $d$ located at location $x_\ell$ at time $t\Delta t$,
- $\hat{h}_\ell^{d,t+1}$ : the distance between vehicles $N_\ell^{d,t}$ and $N_{\ell+1}^{d,t}$ at time $(t+1)\Delta t$,
- $\rho_\ell^{d,t}$ : the mean density in cell $(\ell) \stackrel{\text{def}}{=} [x_\ell, x_{\ell+1}]$ at time t $\Delta t$ :
$$\rho_\ell^{d,t} = \tfrac{N_\ell^{d,t} - N_{\ell+1}^{d,t}}{h}$$

N. Farhi, H. Haj-Salem, M.M. Khoshyaran, J.P. Lebacque, F. Salvarani, B. Schnetzler, F. de Vuyst

- $\hat{\rho}_\ell^{d,t}$ : the mean density of vehicles located between vehicles $N_\ell^{d,t}$ and $N_{\ell+1}^{d,t}$ at time $(t+1)\Delta t$ :

$$\hat{\rho}_\ell^{d,t} = \frac{N_\ell^{d,t} - N_{\ell+1}^{d,t}}{\hat{h}_\ell^{d,t+1}} = \rho_\ell^{d,t} \frac{h}{\hat{h}_\ell^{d,t+1}}$$

The scheme is very simple and illustrated by the following figure

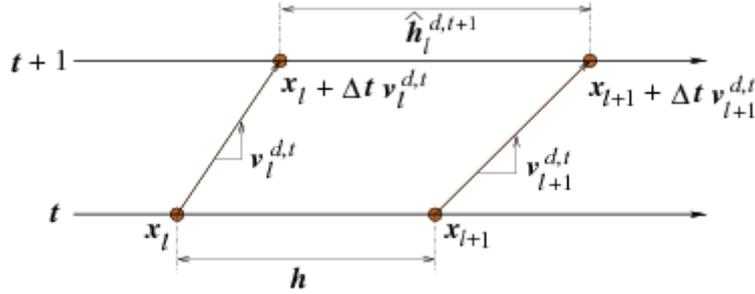

FIGURE 4: Euler-lagrange Remap scheme

We consider a cell $(\ell) = [x_\ell, x_{\ell+1}]$. The number of vehicles (d) which have exited the cell $(\ell)$ and entered the cell $(\ell+1)$ is given by:

$$v_{\ell+1}^{d,t} \Delta t \, \hat{\rho}_\ell^{d,t+1}$$

because the density of these vehicles is $\hat{\rho}_\ell^{d,t+1}$ and the distance covered is $v_{\ell+1}^{d,t} \Delta t$. The underlying hypothesis being that the density is uniform between vehicles $N_\ell^{d,t}$ and $N_{\ell+1}^{d,t}$.

By the same token the number of vehicles (d) which have entered the cell $(\ell)$ and exited the cell $(\ell-1)$ is given by:

$$v_\ell^{d,t} \Delta t \, \hat{\rho}_{\ell-1}^{d,t+1}$$

The conservation of vehicles (d) in cell $(\ell)$ writes:

$$\rho_\ell^{d,t+1} = \rho_\ell^{d,t} + \frac{\Delta t}{h} \left( v_{\ell+1}^{d,t} \hat{\rho}_\ell^{d,t+1} - v_\ell^{d,t} \hat{\rho}_{\ell-1}^{d,t+1} \right) \qquad (8)$$

All that is left is to calculate the $v_\ell^{d,t}$. These are obtained by solving (1) or (2) with respect to the partial densities $\rho_{i,\ell}^{d,t} = \mathfrak{R}_i^d(\bar{\rho}_\ell^t)$. The lane densities are given by (3):

$$\rho_{i,\ell}^t = \sum_{\delta \in D_i} \rho_{i,\ell}^{\delta,t} = \mathfrak{R}_i(\bar{\rho}_\ell^t).$$

The velocities $v_\ell^{d,t}$ are deduced from (4):

$$v_\ell^{d,t} = \sum_{i \in I^d} V_i(\rho_{i,\ell}^t) \rho_{i,\ell}^{d,t} / \rho_\ell^{d,t} \qquad (9)$$

In some sense the approach can be considered lagrangian. Indeed the velocity $v_\ell^{d,t}$ of the vehicles $N_\ell^{d,t}$ is obtained by considering the headway $x_{\ell-1} - x_\ell$.

The scheme is conservative.

The main advantage of the Euler-lagrange remap scheme over, for instance, the Lax-Friedrichs scheme, is the following: the Euler-lagrange remap requires much less flux evaluations

N. Farhi, H. Haj-Salem, M.M. Khoshyaran, J.P. Lebacque, F. Salvarani, B. Schnetzler, F. de Vuyst

1 **Lagrangian Scheme**

2 The conservation law of cars by user classes can be written as follows in lagrangian
3 coordinates.

$$\partial_t r^d + \partial_n v^d = 0, \quad \forall d \in D, \tag{10}$$

5 where $r^d$ and $v^d$ denote respectively the inter-vehicular space and the velocity of class $d$
6 users, and $\partial_t$ and $\partial_n$ denote the partial derivatives with respect to time $t$ and vehicle (or
7 vehicle group) index $n$ respectively.

8 The idea of the scheme is to combine the discretization of the conservation law (10) with a
9 fixed point equation that calculates the split coefficients of the car-densities per user class on
10 different lanes. Let us use the notations:

11 - $\Delta N_n^d$: the size of the group $(n, d)$.
12 - $x_n^d(t)$: the position of the rear of the group $(n, d)$ at time $t\Delta t$.
13 - $\Delta x_{i,n}^d(t)$: the spacing headway of the group $(n, d)$ with respect to its predecessor
14   $(m, e)$ on lane $i$ at time $t\Delta t$. That is $\Delta x_{i,n}^d(t) = x_m^e(t) - x_n^d(t)$.
15 - $\varphi_{i,n}^d(t)$: the split coefficient of the density of group $(n, d)$ over lane $i$, at time $t\Delta t$.

16 The scheme is then given as follows. At each time step $t\Delta t$:

17 1. Determine the split coefficients $\varphi_{i,n}^d(t)$
18

19 $$\varphi_{i,n}^d(t) = \frac{\exp\left[\left\{V_j\left(\frac{\Delta N_n^d \varphi_{i,n}^d(t)}{\Delta x_{i,n}^d(t)}\right) + \theta_i^d\right\}/\nu\right]}{\sum_{j \in I^j} \exp\left[\left\{V_j\left(\frac{\Delta N_n^d \varphi_{j,n}^d(t)}{\Delta x_{j,n}^d(t)}\right) + \theta_j^d\right\}/\nu\right]}, \quad \forall d \in D, \forall i \in I^d. \tag{11}$$

20
21 2. Determine the velocities and the new positions of the groups
22
23 - Estimate the velocities $v_{i,n}^d$ on each lane

$$v_{i,n}^d(t) = V_i\left(\frac{\Delta N_n^d \varphi_{i,n}^d(t)}{\Delta x_{i,n}^d(t)}\right), \quad \forall d \in D, \forall i \in I^d.$$

24 - Take the velocity $v_n^d$ of group $(n, d)$ as follows.

$$v_n^d(t) = \sum_{i \in I^d} \varphi_{i,n}^d(t) \, v_{i,n}^d(t), \quad \forall d \in D.$$

25 - Update the positions of the groups
$$x_n^d(t+1) = x_n^d(t) + \Delta t \, v_n^d(t).$$
26 3. Re-order the groups if necessary.

27

28 **NUMERICAL RESULTS**

29 **Description of the example.**

30 We illustrate in this section, the different schemes and compare them by means of an
31 academic example.

N. Farhi, H. Haj-Salem, M.M. Khoshyaran, J.P. Lebacque, F. Salvarani, B. Schnetzler, F. de Vuyst

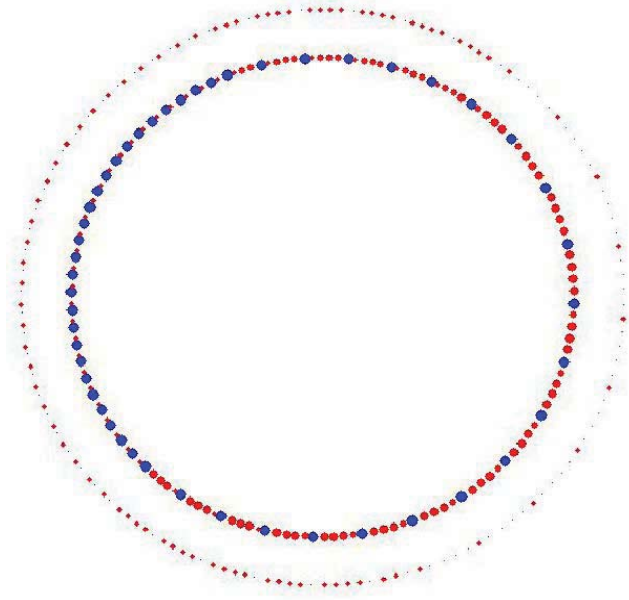

FIGURE 5. Position of vehicle groups on the road. Interior ring: lane 1. The size of a dot corresponds to the proportion of the vehicle group on each lane.

We consider a number of cars moving on a ring road of two lanes ($I = \{1,2\}$). There are two classes of cars ($D = \{1,2\}$). We suppose that the class 1 users move only on lane 1, that is to assume that $\theta_1^1 = +\infty$ and $\theta_2^1 = 0$. Thus we have $\rho_1^1(t) = \rho^1(t)$ and $\rho_2^1(t) = 0$, for all $t \geq 0$. Moreover, we assume that the class 2 users have the same preference for both lanes 1 and 2. More precisely we put

$$\theta \stackrel{\text{def}}{=} \nu^{-1} = 1/12.5 \text{ hr/km}, \text{ and } \quad \theta_1^2 = \theta_2^2.$$

We then have $\varphi_{1,n}^1(t) = 1$ and $\varphi_{1,n}^1(t) = 0$, for all $n, t$. We also have $\varphi_{2,n}^2(t) = 1 - \varphi_{1,n}^2(t), \forall n, t$. In order to simplify the notations, we put $\varphi_n(t) = \varphi_{1,n}^2(t)$. Then $\varphi_{2,n}^2(t) = 1 - \varphi_n(t)$. The fixed point equation is thus written

$$\varphi_n(t) = \frac{\exp\left[\theta\, V_1\left(\frac{\Delta N_n^2\, \varphi_n(t)}{\Delta x_{1,n}^2(t)}\right)\right]}{\exp\left[\theta\, V_1\left(\frac{\Delta N_n^2\, \varphi(t)}{\Delta x_{1,n}^2(t)}\right)\right] + \exp\left[\theta\, V_2\left(\frac{\Delta N_n^2\, (1-\varphi_n(t))}{\Delta x_{2,n}^2(t)}\right)\right]} \qquad (12)$$

which is a fixed point equation in one dimension, easily solved by a standard Newton method.

Below we show some numerical results obtained by applying the Lagrangian scheme.

N. Farhi, H. Haj-Salem, M.M. Khoshyaran, J.P. Lebacque, F. Salvarani, B. Schnetzler, F. de Vuyst

1 **Initial conditions**

2 The initial conditions are those of a Riemann problem. For simplicity, boundary conditions
3 are skipped by considering a circular geometry.
4 The length of the ring is 10 km.
5 The initial conditions are described as a couple of values, left and right hand side of the
6 origin:
7   • ρ1 : [10, 10] veh/km ; ρ 2 : [5, 90] veh/km ; ρ : [15, 100] veh/km
8
9 The figure 6 below describes this initial data.

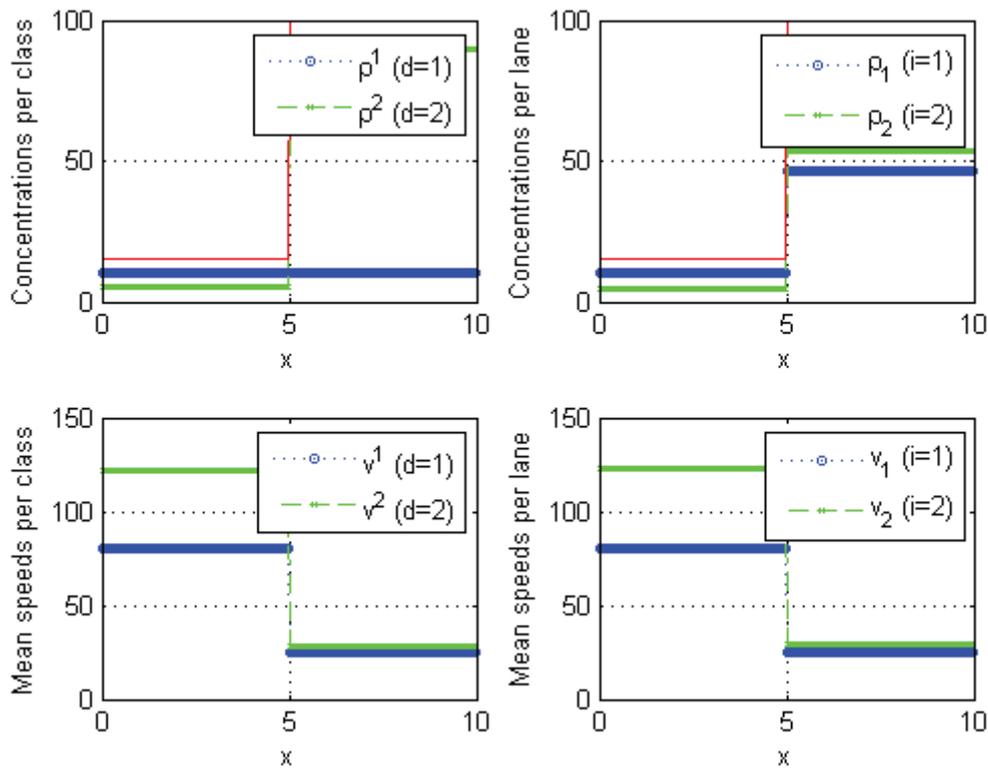

12 FIGURE 6 : Initial conditions of the Riemann problem.
13 Thickest (blue): destination 1, thick (green): destination 2
14

**Results:**

*Lax-Friedrichs scheme.*

The algorithm was used with 1600 cells and a CFL of 1, in order to achieve a good convergence. Nevertheless due to the computational demands of the algorithm it has been deemed more convenient to use in some cases a more efficient version of Lax-Friedrichs, the Rusanov scheme. The Rusanov scheme achieves convergence with 800 cells and a CFL = 0.5 (here we denote CFL the ratio $V_{max}\, \Delta t / \Delta x$).

*Euler-lagrange remap scheme.*

This algorithm requires only 400 cells to achieve an excellent convergence, but a CFL = 0.25. Overall the algorithm is advantageous because the computational demands of each time step are less than those of the Rusanov version of Lax-Friedrichs.
Results are depicted below (Figure 7, snapshot 2 minutes after beginning).

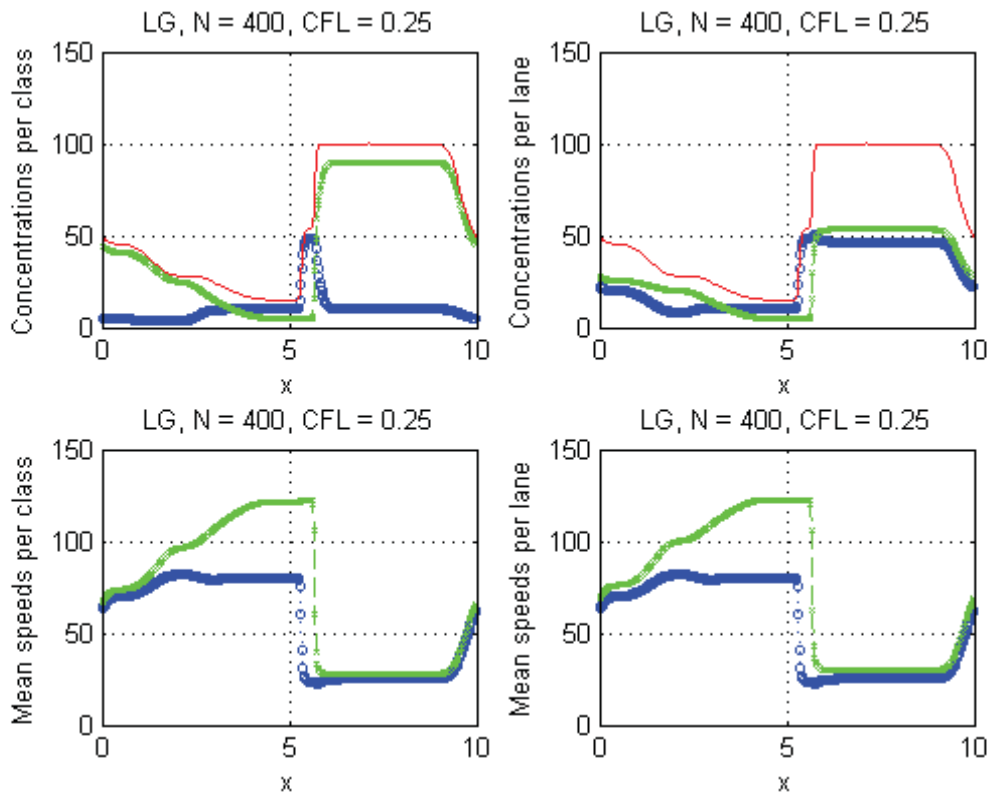

FIGURE 7: Solution of the Riemann problem: densities and speeds per destination class. Rusanov (Lax-Fridriechs) and Euler-lagrange remap
Thickest (blue): destination 1, thick (green): destination 2, thin (red): total (in the case of densities)

N. Farhi, H. Haj-Salem, M.M. Khoshyaran, J.P. Lebacque, F. Salvarani, B. Schnetzler, F. de Vuyst

*Lagrange scheme*

The traffic is discretized into some 146 groups containing 5 vehicles each. Simulation results (after 2 minutes simulated time) are depicted on the figure below.

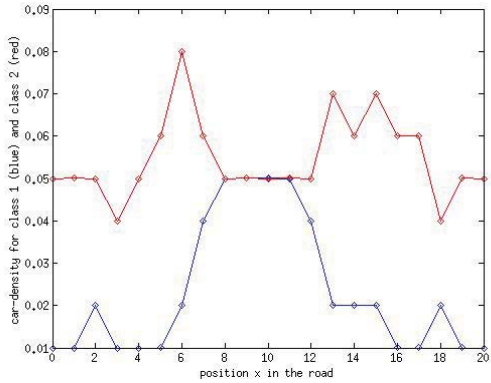
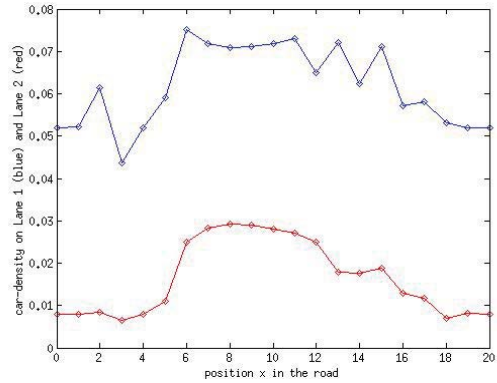

Densities per class　　　　　　　　　　Densities per lane

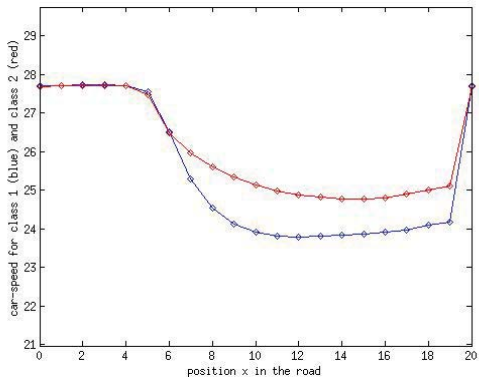
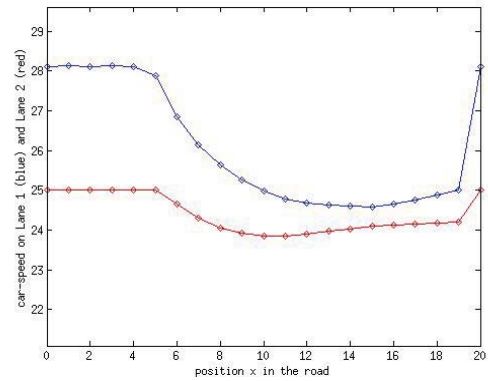

Speeds per class　　　　　　　　　　Speeds per lane

FIGURE 8 : Solution of the Riemann problem: densities and speeds per destination class.
Lagrange scheme
In blue: destination class 1, in red, destination class 2

It is to be noted that the lagrangian scheme is much less precise than the Euler-lagrange remap. This is only to be expected, considering the discretization (much less refined), and some estimates which include significant approximations (for instance the speed and spacing estimates). The potential advantages of lagrangian discretization lie with a correct treatment of intersections (expanding for instance on Khoshyaran Lebacque 2008), and possible extensions to models which accommodate multiple types of vehicles.

**CONCLUSION AND NEXT STEPS**

The analysis of multi-lane data on two data sets has shown that speed is a relevant variable for lane choice and that the traffic measurements are compatible with the Logit lane assignment. The analysis is limited and needs to be pursued, in order to complete the assessment of the

N. Farhi, H. Haj-Salem, M.M. Khoshyaran, J.P. Lebacque, F. Salvarani, B. Schnetzler, F. de Vuyst

validity of the hypotheses underlying the Logit lane assignment model The development of computational tools for the numerical treatment of the model will make this analysis possible at a macroscopic level.

The Euler-Lagrange remap and the Lagrange perform satisfactorily but require improvements. The Euler-Lagrange remap can be improved by considering several lagrangian steps before projecting the solution on the fixed grid. The Lagrange scheme could be modified into that direction too, and also with respect to spacing and speed estimation.

A further step will consist in the extension of the schemes to merges and diverges, following Lebacque Khoshyaran 2008. In the homogeneous case the model simplifies and becomes a model of the GSOM family (see Lebacque et al 2007 for basic properties of this family of models). The new schemes should interface seamlessly with the GSOM model.

In conjunction with this development we expect to investigate further the experimental aspects, in order to improve the physical aspects of the model.

Finally the resulting schemes should complemented with an efficient parametric estimation method and should be used to model an infrastructure such as Boulevard Périphérique, or on the NGSIM data, in order to validate the Logit lane assignment model.

## REFERENCES

Cohen S., (2001), *Aménagement du tronc commun A3-A86 en Seine-Saint-Denis : Evaluation d'impact sur le niveau de service*, Rapport de convention DDE93-INRETS, juin 2001.

Daganzo C.F. (1997) A continuum theory of traffic dynamics for motorways with special lanes. *Transportation Research* B, 31, pp 83-102.

Greenberg J.M., Klar A., Rascle M. (2003) Congestion on Multilane highways. *SIAM Journal of Applied Mathematics* 63(3), 813-818.

Godlewski E., Raviart, P-A. (1996). *Numerical approximation of hyperbolic systems of conservation laws*. Springer Verlag.

Helbing D. (1997). *Verkehrsdynamik*. Springer Verlag.

Holden, H., Risebro, N.H. (2002) *Front tracking for hyperbolic conservation laws*. Springer Verlag.

Hoogendoorn, S.P., Bovy, P.H.L. (1999) Multiclass macroscopic traffic modelling : a multilane generalization using gas-kinetic theory. *Proceedings of the 14th ISTTT (International Symposium on Transportation and Traffic Flow Theory)*, A. Ceder Ed. pp 27-50. Jerusalem.

Jin, W.L. (2006). A kinematic wave theory of lane-changing vehicular traffic. *eprint arXiv:math/0503036*.

Khoshyaran M.M., Lebacque J.P. (2008) Lagrangian modelling of intersections for the GSOM generic macroscopic traffic flow model. *AATT, Athens 2008*.

Khoshyaran M.M., Lebacque J.P. (2012). Numerical solutions to the Logit lane assignment model. To be published in *Procedia Social and Behavioral Sciences*. EWGT 2012, Marne-la-Vallée, France.

N. Farhi, H. Haj-Salem, M.M. Khoshyaran, J.P. Lebacque, F. Salvarani, B. Schnetzler, F. de Vuyst

1   Laval, J., Daganzo, (2005) C.F. Multi-lane hybrid traffic flow model : a theory on the impacts
2   of lane changing maneuvres. *TRB Annual Meeting 2005*.

3   Laval, J., Daganzo, C.F. (2006) Lane-changing in traffic streams. *Transportation Research*
4   40B, 3, pp 251-264.

5   Laval J A and L Leclercq (2008). Microscopic modeling of the relaxation phenomenon using
6   a macroscopic lane-changing model. *Transportation Research Part B*, 42 (6):511-522.

7   Lebacque, J.P., Khoshyaran, M.M. (2002) Macroscopic flow models. Presented at the *6th*
8   *Meeting of the EURO Working Group on Transportation* 1998. Published in "*Transportation*
9   *planning: the state of the art*" (Editors: M. Patriksson, M. Labbé), 119-139. Kluwer Academic
10  Press.

11  Lebacque, J.P., Khoshyaran,M.M. (2009) A stochastic lane assignment scheme for
12  macroscopic multilane traffic flow modelling. *TRISTAN VII*, Tromsø.

13  Lebacque J.P., Mammar S., Haj-Salem H. (2007) Generic second order traffic flow modeling.
14  *Transportation and Traffic Flow Theory 2007* (R.E. Allsop, M.G.H. Bell, B.G. Heydecker
15  Eds). London.

16  Lighthill, M.H., Whitham, G.B. On kinematic waves II : A theory of traffic flow on long
17  crowded roads. *Proc. Royal Soc.* (Lond.) A 229 : 317-345. 1955.

18  Michalopoulos, P.G., Beskos, D.E., Yamauchi, Y. (1984) Multilane traffic flow models :
19  some macroscopic considerations. *Transportation B*, 18, 377-395.

20  Ngoduy, D. (2005) *Macroscopic discontinuity modelling for multiclass multilane traffic flow*
21  *operations*. PhD Thesis, Delft University of Technology, The Netherlands.

22  Richards (1956), P.I. Shock-waves on the highway. *Opns. Res*. 4 : 42-51.

23  Schnetzler, B., Louis, X. Lebacque, J.P., (2010) *A multilane junction model*. Transportmetrica

24